\documentclass[11pt,reqno,a4,authoryear]{imsart}
\usepackage{a4}
\usepackage{graphicx}
\usepackage{amsmath,amssymb,amsthm,color,bbm,ifthen,graphicx,epsf,enumerate,bbm}
\usepackage[authoryear]{natbib}
\usepackage{xspace}
\usepackage[latin1]{inputenc}

\newcommand{\SMdsp}{f_\star}
\newcommand{\SMdspell}{g_\star}
\newcommand{\dsp}{f}
\newcommand{\nbblock}{g}

\newtheorem{thm}{Theorem}

\newtheorem{cor}[thm]{Corollary}

\newtheorem{lem}[thm]{Lemma}
\newtheorem{defi}{Definition}
\theoremstyle{definition}

\newcommand{\note}[2]{}
\newcommand{\LocalLipshitz}[1]{\mathcal{L}_\star(#1)}
\newcommand{\LocalHolder}[3]{\mathcal{F}_\star(#1,#2,#3)}
\def\eqdef{\triangleq}

\def\rset{\mathbb{R}}
\def\rme{\mathrm{e}}
\def\rmi{\mathrm{i}}
\def\rmd{\mathrm{d}}

\def\PE{\mathbb{E}}
\def\PVar{\mathrm{Var}}
\def\PCov{\mathrm{Cov}}
\def\eqsp{\;}

\def\GPH{\mathrm{}}
\def\MSE{\mathrm{MSE}}
\newcommand{\ie}{\emph{i.e.}}

\newcommand{\wrt}{\emph{w.r.t.}}
\newcounter{hyp}

\begin{document}
\begin{frontmatter}
\title{Log-average periodogram estimator of the memory parameter}
\runtitle{Log-average periodogram estimator}
\begin{aug}
\author{Valderio A. Reisen \ead[label=a1]{valderio@cce.ufes.br}}
\author{Eric Moulines \ead[label=a2]{moulines@tsi.enst.fr}}
\author{Philippe Soulier \ead[label=a3]{Philippe.Soulier@u-Paris10.fr}}
\author{Glaura Franco \ead[label=a4]{glaura@est.ufmg.br }}
\address{Departamento de Estatística - CCE - UFES - Av. Fernando Ferrari , s/n - Vitória - ES - Brazil - CEP: 29065-900. \printead{a1}}
\address{GET/Télécom Paris, 46 rue Barrault, 75634 Paris, Cédex 13. \printead{a2}}
\address{Department of Statistics, Universidade Federal de Minas Gerais, Belo Horizonte, MG, Brazil. \printead{a4}}
\end{aug}

\begin{abstract}
This paper introduces a  semiparametric  regression estimator of
the memory parameter for long-memory time series process. It is
based on the  regression in a neighborhood of the zero-frequency
of the periodogram averaged over epochs. The proposed
estimator is theoretically justified and empirical Monte Carlo
investigation gives evidence that the method is very promising to
estimate the long-memory parameter.
\end{abstract}


\end{frontmatter}
\newpage

\section {Introduction}
Let $\{X_t\}$ be a covariance stationary process with spectral density
\begin{equation}
\label{eq:polenzero} \dsp(\omega) \eqdef | \omega |^{-2 d} \SMdsp( \omega) \eqsp, \omega  \in
[-\pi, \pi] \quad \text{and} \quad  d \in (-1/2,1/2) \eqsp,
\end{equation}
where $\SMdsp(\omega)$ is continuous at zero, $0 < \SMdsp(0) <
\infty$ and  $\int \SMdsp(\omega) \rmd \omega < \infty$. To maintain
generality of the short-run dynamics, we do not impose a specific functional form on $\SMdsp(\omega)$.
Equation \eqref{eq:polenzero}  is referred to as a \emph{semiparametric} model for $\dsp(\omega)$,
specifying its form only near zero frequency. Time series with spectral density satisfying \eqref{eq:polenzero}
 can be observed in many areas of applications;
see for example \cite{beran:1994}, \cite{doukhan:oppenheim:taqqu:2003} and the references therein. The process
is said to have \emph{short-memory} when $d=0$, \emph{long-memory} when $d \in (0,1/2)$ and \emph{negative} (or
\emph{intermediate}) memory if $d \in (-1/2,0)$.
Equation \ref{eq:polenzero} is satisfied leading models for long
and negative memory such as Fractional Autoregressive Integrated
Moving Averages (FARIMA) and fractional noise. These, however, are
parametric models, specifying $\dsp (\omega)$ up to finitely many
unknown parameters over all frequencies $(-\pi,\pi]$. The memory
parameter $d$ (like the scale parameter $\SMdsp(0)$) is typically
unknown and should be estimated.  Many works have been devoted to
the estimation of the memory parameter in the semiparametric
context. One of the most popular is the ordinary least squares
(OLS) estimator, due  to \cite{geweke:porter-hudak:1983} (GPH).
This method employs the periodogram  to obtain, through a
regression equation, the OLS estimator of the memory parameter.
Since the GPH estimator has been introduced, other variants of
this method have been suggested with the aim of improving the
quality of estimates and to achieve better asymptotic sample
properties; see \cite{doukhan:oppenheim:taqqu:2003} for an in
depth survey and \cite{nielsen:frederiksen:2005} for a detailed
experimental study.

In some areas of application,  extremely long time series (hundreds of thousands to millions of samples are not
uncommon, for example, in the analysis of teletrafic time-series or in high-frequency finance) has to be dealt
with. Often, the analysis of these data has to be done on-line. Instead of computing the periodogram on the
whole data set, a simple solution to the problem consists in partitioning the sample into subsets, referred to
as \emph{epochs}, computing the periodogram over each epoch and averaging these to obtain an averaged version of
the periodogram. This line of research has been pursued in the parametric context  by \cite{beran:terrin:1994}
but has not, up to our best knowledge, been explored in the semiparametric context.

In this contribution, we study the  averaged periodogram spectral estimator, based  on the division
of series into epochs, to obtain the memory
parameter estimate of a long-memory process. The estimation method
follows the GPH procedure, where the periodogram is replaced by the
averaged periodogram in the regression equation. Some desirable
asymptotic properties of the proposed estimator are derived and
empirical investigation gives  evidence to support the use of
the procedure as an alternative method to GPH to reduce the variance
of the fractional memory parameter.

All these  topics are presented in the paper as follows: Section \ref{sec:model} outlines some properties of the
averaged periodogram. Section \ref{sec:estimation:method} presents the proposed estimator of the memory parameter $d$ and discuss the
statistical performance of the estimator. Section \ref{sec:Monte-Carlo} presents the set-up of the Monte-Carlo experiments and assesses the finite sample properties of our estimator.

\paragraph{Notations:} In the paper,  $a \eqdef b$  denotes that $a$ is defined as $b$, $\lfloor a \rfloor$
is the integer part of $a$ and $a \wedge b \eqdef   \min(a,b)$.

\section{Some properties of the averaged periodogram of a fractional difference process}
\label{sec:model}
 Let  $\{X_t \}_{t=1}^{N}$ be a  realization of the process $X_t$ and  $N =
\nbblock n$, where  $\nbblock$ and $n$ are integer values that correspond , respectively, to the number of
epochs and the sample size of each epoch. The most natural tool for inference in the spectral domain is the
periodogram, defined as the square modulus of the discrete Fourier transform (DFT). Instead of computing the
periodogram over the whole dataset,  the DFT and the periodogram of each individual epoch are obtained at the
Fourier frequencies $\omega_k= 2 \pi k/n$, $1 \leq k \leq [n/2]$, \ie\ for $\ell=0,\dots,\nbblock-1$,
\begin{align}
\label{eq:definition-discrete-Fourier-transform-i}
& d_{\ell,n}(\omega_k) \eqdef (2\pi n)^{-1/2}  \sum_{t=1}^{n} X_{t+ \ell n} \rme^{\rmi t \omega_k} \eqsp, \\
\label{eq:definition-periodogramme-i}
& I_{\ell,n}(\omega_k) \eqdef |d_{\ell,n} (\omega_k)|^2 = (2\pi n)^{-1} \left|
 \sum_{t=1}^{n} X_{t+ \ell n} \rme^{\rmi t \omega_k} \right|^2
 \eqsp.
\end{align}
The averaged periodogram is then defined as the sample mean of the periodogram obtained on the successive epochs,
\begin{equation}
\label{eq:AveragedPeriodogram}
\bar{I}_{\nbblock,n} (\omega_{k}) = g^{-1} \sum_{\ell=0}^{g-1} I_{\ell,n}(\omega_{k})\eqsp.
\end{equation}
The averaged periodogram serves as the basis of the  Welch estimator of the spectral density. Since in the
sequel the spectral density $\dsp$ can be either zero $(d < 0)$ or infinite $(d > 0)$ at zero frequency, we have
found more appropriate to state  results using the normalized periodogram, \ie\ the raw periodogram normalized
by the inverse of the spectral density. For short-memory processes, \ie\  $d=0$ in  \eqref{eq:polenzero}, it is
well known that
\begin{enumerate}[(i)]
\item the periodogram in each epoch is an asymptotically unbiased estimate of the spectral
  density, \ie\ $ \dsp^{-1}(\omega_k)  \PE[I_{0,n}(\omega_{k})] = 1 + O(n^{-1})$, $1 \leq k \leq [n/2]$,
   where the $O(n^{-1})$ term is uniform with respect to ( \wrt\ )
  the frequency coordinates $k$,
\item the DFT ordinates in a single epoch at distinct Fourier frequencies are asymptotically uncorrelated,
\begin{multline*}
   \PCov \left( \dsp^{-1/2}(\omega_j) d_{0,n}(\omega_{j}), \dsp^{-1/2}(\omega_k) d_{0,n}(\omega_{k}) \right) = \delta_{j,k}  + O(n^{-1}) \eqsp, \\ 1 \leq j \leq k < [n/2],
  \end{multline*}
 where the $O(n^{-1})$ term is uniform \wrt\ the frequency indices $j,k$.
\item the DFT ordinates over different epochs are asymptotically uncorrelated, i.e. for any $\ell > 0$,
\begin{multline*}
   \PCov \left( \dsp^{-1/2}(\omega_j) d_{0,n}(\omega_{j}), \dsp^{-1/2}(\omega_k) d_{\ell,n}(\omega_{k}) \right) =  O(n^{-1}) \eqsp, \\ 1 \leq j \leq k < [n/2],
  \end{multline*}
where the $O(n^{-1})$ term is uniform \wrt\ the frequency indices $j,k$ and the epoch index $\ell$.
\end{enumerate}
Assuming that the number of epochs $\nbblock$ is a fixed integer, the above
results  imply that
\begin{itemize}
\item[(iv)]  the averaged periodogram is an asymptotically unbiased estimator of the spectral density, \ie\
$\dsp^{-1}(\omega_k) \PE[\bar{I}_{\nbblock,n}(\omega_k)] =  1 + O(n^{-1})$, where the $O(n^{-1})$ term is uniform in $k$;
\item[(v)] the averaged periodogram ordinates are asymptotically uncorrelated and its variance is equal
 to the square of the spectral density divided by the number of epochs,
\[
\PCov\left(\dsp^{-1}(\omega_j) \bar{I}_{\nbblock,n}(\omega_j), \dsp^{-1}(\omega_k) \bar{I}_{\nbblock,n}
(\omega_{k})\right) = g^{-1} \delta_{j,k}  + O(n^{-1}).
\]
\end{itemize}

We now discuss  some  periodogram properties for the process  \eqref{eq:polenzero} with $d \neq 0$.

It was first shown by \cite{kunsch:1986} and  exhaustively investigated by \cite{hurvich:beltrao:1993} in the
single epoch case that the asymptotic behavior of low frequency DFT ordinates departs strongly from the weak
dependence situation: in particular, the DFT ordinates computed at Fourier frequencies $\omega_j = 2 \pi j/n$
and $\omega_k = 2 \pi k /n$ for fixed positive integers $j$ and $k$ are  correlated as $n \to \infty$. Not
surprisingly, the same departure from the weak dependence behavior is also observed in  a multiple epochs
scenario. The dependence in particular implies that the correlation among DFT coefficients over different epochs
does not asymptotically vanish. More precisely, the following theorem shows that the correlation of the DFT coefficients evaluated at Fourier frequencies
$\omega_j$ and $\omega_k$, for \emph{fixed} $j$ and $k$, computed over different epochs does not vanish as $n \to \infty$.
\begin{thm}
\label{theo:covariance-DFT-multiple-epochs}
Assume that the spectral density $\dsp$ is given by \eqref{eq:polenzero}. \note{EM}{there is a missing condition there}
Let $\ell \geq 0$ and $1 \leq j \leq k$ be fixed positive integers. Then,
\begin{align}
&\lim_{n \to \infty} \omega_j^d \omega_k^d \PE \left[ d_{0,n}(\omega_j) \bar{d}_{\ell,n}(\omega_k) \right]
= \frac{(2 \pi j)^d (2 \pi k)^d}{2 \pi}  D_1(d;j,k, \ell) \eqsp \\
&\lim_{n \to \infty} \omega_j^d \omega_k^d \PE \left[d_{0,n}(\omega_j) d_{\ell,n}(\omega_k) \right]= \frac{(2
\pi j)^d (2 \pi k)^d}{2 \pi}  D_2(d;j,k, \ell) \eqsp.
\end{align}
where
\begin{align}
&D_1(d;j,k,\ell) \eqdef \int_{-\infty}^\infty |\omega|^{-2d} \Delta(\omega - 2 \pi j) \Delta (2 \pi k - \omega) \rme^{-\rmi \ell \omega} \rmd \omega \eqsp, \\
&D_2(d;j,k,\ell) \eqdef  \int_{-\infty}^\infty |\omega|^{-2d} \Delta(\omega - 2 \pi j) \Delta (-2 \pi k - \omega) \rme^{-\rmi \ell \omega} \rmd \omega \eqsp
\end{align}
with $\Delta(\omega) \eqdef (\rme^{\rmi \omega} - 1)/(\rmi \omega)$.
\end{thm}
\begin{proof} See Section \ref{proof:theo:covariance-DFT-multiple-epochs}. \end{proof}

Theorem \ref{theo:covariance-DFT-multiple-epochs} combined with the results of Theorem 5 given in
\cite{hurvich:beltrao:1993} implies that the averaged periodogram is an asymptotically  biased estimate of the
spectral density, \ie\ $\lim_{n \to \infty} \dsp^{-1}(\omega_j) \bar{I}_{g,n}(\omega_j) \ne 1$ for any given
 $j$,  $1\leq j \leq k$, and  the correlation of the averaged periodogram ordinates does not asymptotically vanish,
\ie\
$$\lim_{n \to \infty} \dsp^{-1}(\omega_j) \dsp^{-1}(\omega_k)
\PCov(\bar{I}_{g,n}(\omega_j),\bar{I}_{g,n}(\omega_k)) \ne \delta_{j,k} \eqsp, $$
for any given $j$ and $k$.
Nevertheless, under appropriate regularity condition for the spectral density of the short-memory process
$\SMdsp$, it can established that there exist two sequences $\{ r(\dsp;k) \}$ and $\{ r(\dsp;j,k) \}$
satisfying,
 for all $n$ and all $1 \leq j < k \leq [n/2]$,
\begin{align*}
&| \PE[\bar{I}_{g,n}(\omega_k)/\dsp(\omega_k)] - 1| \leq r(\dsp;k) \eqsp, \\
&|\PCov( \dsp^{-1}(\omega_j) \bar{I}_{g,n}(\omega_j), \dsp^{-1}(\omega_k) \bar{I}_{g,n}(\omega_k) | \leq r(\dsp;g,k,j) \eqsp,
\end{align*}
such that
\begin{itemize}
\item $\lim_{k \rightarrow \infty} r(\dsp;k) = 0$, which means that the bias is small for frequencies
sufficiently far away from  zero, and
\item for any sequence $m_n$ such that $\lim_{n \to \infty} m_n= \infty$ and $\lim_{n \to \infty} m_n/n < 1/2$,
$\sum_{1 \leq j < k \leq m_n} r(\dsp;g,j,k) = O\left(\log^r(m_n) \right)$, for some $r > 0$. This shows that
 whereas the dependence does not vanish it is small when the Fourier frequencies are sufficiently
far apart.It is also possible (see the results below) to get bounds of this quantity as a function of $g$.
\end{itemize}
To derive the above results,  some regularity assumptions  need to be imposed on the spectral density of the
short-memory process. Consider the  set of functions, which is adapted from
\cite[Theorem 2]{robinson:1995} ( see also \citep{soulier:2001} and \citep{moulines:soulier:2003}).
\begin{defi}
\label{eq:definition-L}
For $\mu \geq 1$ and $\nu \in (0,1]$. Let $\LocalLipshitz{\mu,\nu}$ be the set of functions $\phi: [-\pi,\pi] \to \rset_+$ satisfying
for all $\omega,\omega' \in [-\pi,\pi] \setminus \{0\}$
\begin{align}
&\max_{\omega \in [-\pi,\pi]} \phi(\omega) \leq \mu \phi(0)\eqsp,\\
\label{eq:Lipshitz-Condition-0}
&|\phi(\omega)- \phi(\omega')| \leq \mu  \phi(0) \frac{ |\omega-\omega'|}{|\omega| \wedge |\omega'|} \eqsp, \\
\label{eq:Lipshitz-Condition-1}
&|\phi(\omega) - \phi(\omega')| \leq \mu \phi(0) \left||\omega|-|\omega'|\right|^\nu \eqsp,
\end{align}
\end{defi}
The set of functions $\LocalLipshitz{\mu,\nu}$ contains all the functions which are strictly positive and continuously differentiable on $[-\pi,\pi]$. In particular, if $\SMdsp$ is the spectral density of an stationary and invertible ARMA process, then $\SMdsp \in \LocalLipshitz{\mu,1}$. More interestingly, it also contains functions of the form $\tilde{\SMdsp}(\omega) = \SMdsp(\omega) + \sigma^2 |\omega|^\nu$, where $\SMdsp$ is
 a strictly positive continuously differentiable function. Note that such spectral density appears in the so-called signal-plus-noise model,
 where a fractional process with smooth spectral density $\SMdsp$ is observed in presence of white noise, uncorrelated from the process.
\begin{thm}
\label{theo:covariance-DFT-multiple-epochs-bounds}
Let $\delta_-,\delta_+ \in [0,1/2)$, $\mu \geq 1$ ,and $\nu \in (0,1]$ be constants. Assume that $\dsp$ is given by \eqref{eq:polenzero}
with $\SMdsp \in  \LocalLipshitz{\mu}$ and $d \in [-\delta_-,\delta_+ ]$.
Then, there exists a constant $C$ (depending on the constants $\delta_-$, $\delta_+$ and $\mu$) such that, for any $1 \leq j  < k \leq n/4$ and $\ell \geq 0$,
\begin{equation}
\label{eq:biais-periodogram} \left| \PE[ \dsp^{-1}(\omega_k) I_{0,n}(\omega_k)] - 1 \right| \leq C \log(1+k)/k
\eqsp,
\end{equation}
\begin{multline}
\label{eq:covariance-DFT} \dsp^{-1/2}(\omega_j) \dsp^{-1/2}(\omega_k) \left|\PE [ d_{0,n}(\omega_j) \bar{d}_{\ell,n}(\omega_k)] +
\PE[d_{0,n}(\omega_j) d_{\ell,n}(\omega_k)] \right|  \\
\leq  C \log(k) j^{-|d|}k^{|d|-1} \left( \ell^{-\{ (1-2d) \wedge 1\}} + (\ell n)^{-\nu} \right)\eqsp.
\end{multline}
\end{thm}
\begin{proof}
See section \ref{sec:theo:covariance-DFT-multiple-epochs-bounds}.
\end{proof}
It is worthwhile to note that the dependence among  successive epochs does not asymptotically vanish as $n \to
\infty$, in strong contrast with the short-memory case. Also, the strength of the dependence among the epochs
depends on the memory coefficient $d$.

 Using Corollary 2.1 in \cite{soulier:2001} and under the additional
assumption that the process is Gaussian, it is possible to
translate the results above to non-linear transforms of the DFT
ordinates, for instance to the  "$\log$" function of the
average periodogram. These are presented in the following
corollary which lead to  results that provide  a theoretical
justification for the estimator proposed in this paper. Some
additional notations are required to state the results. Let $U$ be
a central chi-square with $2g$ degrees of freedom. Then, $\PE[
\log(U/2)] = \psi(g)$ and $\PVar[ \log(U/2) ]= \psi'(g)$ where
$\psi$ is the digamma function (see \cite[p.~198]{johnson:kotz:1970}). For instance, $\psi(1)= - \gamma$, where
$\gamma$ is the Euler constant and $\psi'(1)= \pi^2/6$. It is well-known that
$\lim_{g \rightarrow \infty} g \psi'(g) = 1 $. Hence, for large $g$, $
\PVar[ \log(U/2)] $ = $O(g^{-1})$. Let
\begin{equation}
\label{eq:decomposition-averaged-periodogram} \xi_{n,k} \eqdef \log\left[ \bar{I}_{\nbblock,n}(\omega_k) \right]
-  \log\left[\SMdsp(\omega_k)\right] - \psi(\nbblock) + \log(\nbblock)  \eqsp.
\end{equation}
The following corollary   establishes the statistical properties of $\xi_{n,k}$.

\begin{cor}
\label{cor:decomposition-averaged-periodogram}
  Assume that $X_t$ is a Gaussian process. Then, there exists an integer $K$ and a constant
$C$ (depending only on the constants $\delta_-$, $\delta_+$, $\mu$, and $K$) such that, for any $K \leq j  < k
\leq [(n-1)/2]$
\begin{align*}
& \left| \PE[ \xi_{n,j}] \right| \leq C \log(1+j)/j \eqsp, \\
&\left| \PVar[\xi_{n,j}] - \psi'(\nbblock) \right| \leq C \log^2(1+j)/j^2 \eqsp,  \\
&  |\PCov( \xi_{n,j}, \xi_{n,k})|  \leq  C \log^2(k) j^{-2|d|} k^{2|d|-2} \eqsp.
\end{align*}
\end{cor}
\begin{proof}
See Section \ref{sec:theo:covariance-DFT-multiple-epochs-bounds}.
\end{proof}

\section{Estimation of the memory parameter based on log-periodogram regression}
\label{sec:estimation:method} As an application of the results obtained above, it is argued in this section that
the averaged periodogram is a simple mean to reduce the variance of semiparametric estimator of the memory
parameter based on log-periodogram regression. For simplicity,  in this contribution the focus will be on the
GPH estimator proposed by \cite{geweke:porter-hudak:1983} (GPH) and further analyzed by \cite{robinson:1995} and
\cite{hurvich:deo:brodsky:1998}. The same reduction of variance holds for the bias reduced log-periodogram
estimator introduced by \cite{andrews:guggenberger:2003}, which is based on regression of $\log \SMdsp(\omega)$
by an even polynomial of degree $2r$,  and for the estimator introduced by \cite{guggenberger:sun:2006}, which
is obtained by taking a weighted average of GPH estimators over different bandwidths.

The GPH estimator is the ordinary least square (OLS) regression estimator obtained from an approximated
regression equation of the logarithm of the spectral density \eqref{eq:polenzero}, having the logarithm of the
spectral density as the dependent variable and $\log ( \omega )$ as the independent
variable. Taking the logarithm of \eqref{eq:polenzero}, the log-spectral density can be expressed as
\begin{equation}
\log \dsp(\omega) = \log \SMdsp(0) - 2 d \log (\omega) +
\log \left[ \SMdsp(\omega) / \SMdsp(0)  \right] \eqsp.
\end{equation}
For the proposed estimator, the spectral density $\dsp(\omega)$
 is replaced by  the averaged periodogram $\bar{I}_{\nbblock,n}(\omega)$,
and using the decomposition \eqref{eq:decomposition-averaged-periodogram},
an  estimate of $d$ is obtained from the  regression equation
\begin{equation}
\label{eq:emc}
\log \left[ \bar{I}_{\nbblock,n}(\omega_k) \right]= a_{0} -2d \log (\omega_k)
 + \log\left[ \SMdsp(\omega_k)/ \SMdsp(0) \right] + \xi_{n,k}\eqsp,
\end{equation}
where the intercept is $a_{0}=  \log \SMdsp(0)  + \psi(\nbblock) - \log(\nbblock)$ and the random variables $\{ \xi_{n,k} \}$
are defined in Corollary \ref{cor:decomposition-averaged-periodogram}.
The GPH estimate of the memory parameter $d$ is thus given by
\begin{equation}
\label{eq:definition-GPH-estimator} \hat{d}^{\GPH}_{m_n,g} =
\sum_{k=1}^{m_n} a_{k,n}(m_n) \log\left[ \bar{I}_{\nbblock,n}
(\omega_k) \right]
\end{equation}
where $\{ m_n \}$, the bandwidth in the regression equation \eqref{eq:emc}, is a sequence of integers and  the weights
$a_{k}(m_n)$ are given by
\begin{equation}
\label{eq:definition-weight-GPH}
a_{k,n}(m_n) \eqdef \frac{[-2 \log(\omega_k)] - m_n^{-1} \sum_{j=1}^{m_n}
[-2 \log (\omega_j)]}{\sum_{k=1}^{m_n} \left\{[-2\log(\omega_k)] - m_n^{-1} \sum_{j=1}^{m_n} [-2 \log(\omega_j)] \right\}^2} \eqsp.
\end{equation}
We will now derive a central limit theorem for the above estimator. To do this, it is required to state some
additional regularity conditions of the spectral density of the short-memory process; see \cite{giraitis:robinson:samarov:2000}.

\begin{thm}
\label{theo:gphgauss}
Assume that $\{ X_t \}$ is a Gaussian process with spectral density $\dsp$ satisfying \eqref{eq:polenzero} with
$\SMdsp \in  \LocalLipshitz{\mu,\beta}$ for some $\mu < \infty$ and $\beta \in (0,1]$ and
\[
\text{for all} \quad \omega \in [-\Omega_0,\Omega_0], \quad |\SMdsp(\omega) - \SMdsp(0)|  \leq \mu \SMdsp(0) |\omega|^{\beta} \eqsp.
\]
Let $\{m_n\}$ be a non-decreasing sequence of integers such that
\begin{equation}
\label{eq:ratemgph}
\lim_{n\rightarrow\infty} ( m_n^{-1} + m_n^{2 \beta+1} n^{-2 \beta} ) = 0 \eqsp.
\end{equation}
Then $\sqrt{m_n} ( \hat{d}^{\GPH}_{m_n,g} - d )$ is asymptotically
 distributed as Gaussian with zero-mean and variance $\psi'(g)/4$.
\end{thm}
\begin{proof}
see section \ref{sec:proof:theo:gphgauss}.
\end{proof}

Using similar arguments  as those  given in \cite{hurvich:deo:brodsky:1998},  the bias and variance of
$\hat{d}^\GPH_{m_n,g}$ are computed by assuming that $\SMdsp$ is three times differentiable in a neighborhood of
the zero frequency. Using once again Corollary \ref{cor:decomposition-averaged-periodogram},
\begin{equation}
\label{eq:bias-GPH} \PE \left[ \hat{d}^{\GPH}_{m_n,\nbblock}\right] - d  = \sum^{m_n}_{j = 1} a_{j,n}(m_n)
\log \SMdsp(\omega_{j})+ \sum^{m_n}_{j = 1} a_{j,n}(m_n) \PE[\xi_{n,j}]
\end{equation}
where $\xi_{n,j}$ is defined in
\eqref{eq:decomposition-averaged-periodogram} and
\begin{equation}
\label{eq:variance-GPH} \PVar( \hat{d}^{\GPH}_{m_n,\nbblock}) =  \sum^{m_n}_{j = 1} a^{2}_{j,n}(m_n) \PVar
(\xi_{n,j}) +  \sum^{m_n}_{k = 1}\sum^{m_n}_{j = k + 1} a_{j,n}(m_n) a_{k,n}(m_n)
\PCov(\xi_{n,j},\xi_{n,k}) \eqsp.
\end{equation}
Along the same lines as \cite{hurvich:deo:brodsky:1998} (Lemma 1 to Lemma 8), we establish an explicit
expression for the mean-square error (MSE) of the proposed estimator.
\begin{thm}
\label{prop:MSE-GPH} Assume that $\SMdsp \in  \LocalLipshitz{\mu}$
 and satisfies the conditions  $\SMdsp^{'}(0) = 0$, $ \left| \SMdsp''(\omega)\right| < \infty $
and $\left|\SMdsp'''(\omega) \right|< \infty$ for any $\omega \in [-\Omega_0, \Omega_0]$ where $\Omega_0 \in (0,\pi]$.
Then,
\begin{equation}
\label{eq:bias-GPH-1} \PE[ \hat{d}^{\GPH}_{m_n,\nbblock} - d ] = -
\frac{2\pi^{2}\SMdsp''(0)m_n^{2} } {9 \SMdsp(0)n^{2}} + o\left(
\frac{m_n^2}{n^2}\right) + O\left( \frac{\log^3(m_n)}{m_n} \right)
\end{equation}
and
\begin{equation}
\label{eq:variance-GPH-1}
\PVar(\hat{d}^\GPH_{m_n,\nbblock}) =  \frac{\psi'(\nbblock)}{4m_n} + o \left( \frac{1}{m_n} \right) \eqsp.
\end{equation}
\end{thm}
Neglecting the remainder terms in the bias and  variance, and assuming that $\SMdsp''(0) \ne 0$ minimizing the
approximate expression for the MSE, i.e.

\begin{equation}
\mathrm{MSE}(n,\nbblock) =
\left[ \frac{2\pi^{2}\SMdsp''(0)m_n^{2} } {9 \SMdsp(0)n^{2}} \right]^2 +
\frac{\psi'(\nbblock)}{4m_n}
\end{equation}
with respect to the bandwidth parameter $m_n$ for a given number
of epochs $\nbblock$ yields the asymptotically optimal choice for
the bandwidth $m_n(\nbblock)$
\begin{equation}
\label{eq:optimal-bandwith} m_n(\nbblock) \eqdef \left( \frac{\psi'(\nbblock)}{16 B_\star} \right)^{1/5}  n^{4/5} \eqsp,
\end{equation}
where  $B_\star \eqdef (4/81) \pi^4 \left\{ \SMdsp''(0)/\SMdsp(0) \right\}^{2}$. With this choice for $m_n(\nbblock)$,
the optimal value for the mean-square error is
\begin{equation}
\label{eq:optimal-value-mean-square}
\mathrm{MSE}(n,\nbblock) = C_\star  \{ \psi'(\nbblock) \}^{4/5} n^{-4/5} \eqsp,
\end{equation}
where $C_\star \eqdef \left\{ (16 )^{-4/5}  B_\star^{-2/5}+ (16
B_\star)^{1/5}/4  \right\}$. We will now discuss the potential
advantages in performance obtained by dividing the series into
epochs. The optimal MSE of the classical GPH (using a single
epoch) is given by $\MSE^\GPH(N,1)= C_\star \{ \psi'(1) \}^{4/5}
N^{-4/5} $.

Dividing  the series into $\nbblock$ epochs each of size $n= N/\nbblock$, the optimal MSE is given by
\begin{equation}
\label{eq:optimal-MSE-g-blocks}
\mathrm{MSE}^\GPH(N/\nbblock,\nbblock)=  C_\star  \{ \psi'(\nbblock) \}^{4/5} (N/\nbblock)^{-4/5}= C_\star \{ \nbblock \psi'(\nbblock) \}^{4/5} N^{-4/5} \eqsp.
\end{equation}

Since $\nbblock \mapsto \nbblock \psi'(\nbblock)$ is a decreasing function, the optimal MSE is also a decreasing
function of the number of epochs, which means that \emph{dividing  the series into epochs is a very simple way
to improve the MSE}. The quantity $\nbblock \psi'(\nbblock)$ decreases from $\pi^2/6$ to $1$ and as $\nbblock$
goes to infinity (more precisely,  $\nbblock \psi'(\nbblock)= 1 + 1/(2 \nbblock) + O( \nbblock^{-2}) $). For
$\nbblock=3$,  $\nbblock \psi'(\nbblock)$ is 1.1848 and its value  changes slowly thereafter, as it is shown in
the following tabulation.

\begin{table}[h]
\centering
\begin{tabular}{|c|c|c|c|c|c|c|c|c|}\hline\hline
$m$ & 1 & 2 & 3 & 4 & 5 & 6 & 7 & 8\\\hline
$m\psi'(m)$ & 1.646 & 1.290 & 1.185 & 1.138 & 1.108 & 1.090 & 1.080 & 1.070 \\\hline \hline
\end{tabular}
\end{table}

\section{Monte-Carlo results}

 \label{sec:Monte-Carlo} This section
provides a limited Monte Carlo experiment to support our claims. For this purpose, realizations of a Gaussian white
noise sequence $\varepsilon _{t}$,  $t =1,\cdots,n$, with unit variance, were generated by IMSL-FORTRAN
subroutine DRNNOR and trajectories of Gaussian processes $\left\{ X_{t}\right\}$ with spectral density
satisfying \eqref{eq:polenzero} were simulated according to the procedure outlined by \cite{hosking:1981}.
 To assess the performance of  $\hat{d}_{m_n,\nbblock}$,
  we compute the bias, the mean-square error ($mse$) and the coverage rates ($cr$) of
  the asymptotic confidence interval
 based on the normal distribution (see Theorem \eqref{theo:gphgauss}). The quantities were calculated
 based on 2,000 replications for different sample sizes $N$ and  number of epochs $\nbblock$. The results
 are displayed in Tables \ref{tab:result1} to \ref{tab:result3}.  In each experiment, the
sampling distribution for the standardized $\hat{d}_{m_n,\nbblock}$ estimator was calculated to obtain the
 coverage rate, which refers to the percentage of cases where the true value of $d$ ($d$ = 0.3) lies inside  the
95\% asymptotic confidence interval ($\hat{d}_{m_n,\nbblock} \pm 1.96 \sigma_{e,n}$), where $\sigma^2_{e,n}$ is the asymptotic
variance. We use two different approximations of the asymptotic variances,
which are asymptotically equal but different for finite sample size.
The variances are  $\sigma_{a,n}^{2}= \psi'(\nbblock)/4m_n$ (which is the limiting variance in
Theorem \ref{theo:gphgauss})
and
\[
\sigma_{r,n}^{2}={\frac{\psi'(\nbblock)}{\sum_{k=1}^{m_n} \left\{[- 2 \log(\omega_k)] - m_n^{-1} \sum_{j=1}^{m_n} [-2 \log(\omega_j)]
\right\}^2}}
\]
which is the variance of the regression obtained as in the case where  the averaged periodogram ordinates are
independent with equal variance $\psi'(\nbblock)$.  For the ARFIMA process with short-memory parameters, the bandwidth
$m_{n}$ is the one given in the previous section ($m_{n}(\nbblock)$) that minimizes the asymptotic mean-square error of the estimator;
we assume here that the parameters specifying the short-memory component $\SMdsp$ of the ARFIMA model are known,
 to avoid discussing the separate issue of the optimal choice of the bandwidth (this may be seen as an oracle
 estimator in this semiparametric context). In the case where $\SMdsp( \omega) \equiv 1$,  short-memory dynamics
  are not present, hence two large and fixed
bandwidths were used, $m_{1}= n^{0.7}$ and $m_{2}= \frac{n-1}{2}$. These bandwidth choices were considered with
the aim to verify the finite property of the bandwidth  on the estimates and  the convergence of the
standardized estimator to the normal distribution.

Results from Table \ref{tab:result1} support the asymptotic properties discussed in the previous sections for the
ARFIMA$(0,d,0$) model. As it can be observed, the mean-square error decreases as the number of epochs  $\nbblock$ increases.
Breaking the series in a fixed number of epochs can produce a significant reduction in the mean-square error (this effect is similar to the pooling
in frequency domain advocated in \cite{robinson:1995}).
Although large $\nbblock$ does not bring too much gain in terms of the mean-square error,
it does not cause penalty in the estimates unless it reaches very large value, as it has been already discussed in the previous section.
%
%

 Even if the sample sizes used in this limited Monte-Carlo experiments are not large,
  the coverage rates of the asymptotic confidence intervals are reasonably accurate.
  The bandwidth $m_{1}$ produces estimates which are (as expected) less accurate  than $m_{2}$.
  For both bandwidths, the coverage of the asymptotic confidence intervals
  based on $\sigma_{r,n}$  are precise, even for relatively small sample sizes.
 The coverage rates of the standardized estimators using the asymptotic standard deviation ($cr_a$)
 is reduced as $\nbblock$ increases. This indicates that the convergence of $\hat{d}_{m_n,\nbblock}$
 , standardized by   $\sigma_{a}$,   to the N(0,1) is affected by the sample size reduction.

\begin{table}[h]
\centering \caption{ARFIMA($0,d,0$)}{\footnotesize
\begin{tabular}{c c c c c cc|| cccccc }
  \hline
  \hline
    &   &   &\multicolumn{7}{c}{$d=0.3$}   \\
    \cline{4-12}
  $N$ & $g$& $m_1$   & $mean$ & $mse$ & $cr_{r}$& $cr_a$ & $m_2$& $mean$ & $mse$ & $cr_{r}$ &$cr_a$\\
  \hline
  512 & 1  &  78  & 0.3032 & 0.00681  & 94.7& 92.1 &  255 & 0.3035 & 0.00222  & 96.0& 97.9\\
      & 2   & 48 &  0.3068  & 0.00443 & 95.7& 91.5 & 127 &  0.3051 & 0.00181    & 95.2& 95.6\\
      & 4   & 29& 0.3073  & 0.00371  & 95.2& 87.4  & 63& 0.3034  & 0.00193  & 94.0&91.1\\
 \hline
 2048 & 1  &  207 & 0.3004& 0.00215  & 95.2 &94.6  &  1023 & 0.3006& 0.00052  & 94.6& 98.6\\
      & 2   & 128&  0.3027 & 0.00148 & 95.8& 92.7& 511 &  0.3007 & 0.00041    & 94.8&98.1\\
      & 4   & 78& 0.3032  & 0.00109  & 96.2 & 93.1& 255& 0.3027  & 0.00036  & 96.0& 97.4\\
      & 8  & 48& 0.3034  & 0.00096  & 94.5  &  89.6& 127& 0.3022  & 0.00039  & 94.8&95.9\\
      & 16  & 29& 0.3043  & 0.00084  & 95.7 &  88.8& 63& 0.3047  & 0.00046 & 94.1&90.7\\
  \hline
 \hline
 8192 & 1  &  548 &0.3007& 0.00078  &94.5&94.3 &  4095 & 0.3007& 0.00013  & 95.2& 99.6\\
      & 2   & 337& 0.3021& 0.00054 & 95.6& 94.4& 2047 &  0.3006 & 0.00010    & 94.9&99.4\\
      & 4   & 207&0.3016& 0.00039 & 94.3 & 93.0& 1023& 0.3008  & 0.00009  & 94.9& 99.0\\
      & 8  &128& 0.3036 & 0.00032  & 94.1 &  92.3& 511& 0.3016  & 0.00009  & 94.3&98.3\\
      & 16  & 78& 0.3034 &0.00027 & 94.2 & 90.3& 255& 0.3021  & 0.00009 & 94.2&97.3\\
  \hline

\end{tabular}}
\label{tab:result1}
\end{table}

Table \ref{tab:result2} gives the results for the ARFIMA process when non-trivial short-memory
components are present in the model.  As it has
often been reported in the long-memory literature,
 the short-memory component  causes significant bias in the
estimator of the long-memory parameter, especially if the bandwidth is not properly set; see, for example,
\cite{hurvich:beltrao:2004}, \cite{hurvich:ray:1995} and \cite{reisen:abraham:lopes2001} among others. As
shown by the asymptotic analysis,  the bias of the estimator is not significantly affected by the division of the sample
into epochs, \ie\ at least when the number of epochs $\nbblock$ is not large compared to $N$.
 This experiment shows that the improvement of the estimates in terms of mean-square error
 depends on the sample size $N$ and the number $\nbblock$ of epochs. For all cases, a decrease of the mean-square error
 is observed when taking $\nbblock=2,3$,  but the improvement becomes marginal when $\nbblock \geq 4$ . This empirical
property of the estimator is not surprising and it was justified theoretically in the previous section.
For example,  for $\nbblock=2$ the decrease of the $mse$ predicted from the asymptotic expression
\eqref{eq:optimal-MSE-g-blocks} when using the optimal bandwidth is $0.783$, which is consistent with the values
found in the Monte-Carlo experiments. The coverage rate is stable when using the regression variance in the
standardized estimator, by the other hand  the values of  $cr_a$   are reduced
as $\nbblock$ increases, which is presumable due   to the sample size reduction.


\begin{table}[h]
\centering \caption{Estimation in ARFIMA($1,d,0$) model using the optimal  bandwidth } {\footnotesize
\begin{tabular}{c c c c c cc||ccccc }
  \hline
  \hline
    &   &\multicolumn{3}{c}{$d=0.3$, $\phi= -0.3$}&&&\multicolumn{3}{c}{$d=0.3$, $\phi =0.3$}& & \\
    \cline{2-12}
  $N$ & $g$& $m_n$   & $mean$ & $mse$ & $cr_{r}$ &  $cr_{a}$ &$m_n$& $mean$ & $mse$ & $cr_{r}$ & $cr_{a}$\\
  \hline
  512 & 1  &  103 & 0.26682 & 0.00574  & 92.8&90.2  &  62 & 0.3440 & 0.01083  & 91.7& 87.3\\
      & 2   & 49 &  0.26973  & 0.00511 & 92.7& 88.6 & 29 &  0.3445 & 0.01049    & 91.8& 83.3\\
      & 4   & 23& 0.27282  & 0.00572 & 92.7&  84.9 & 14& 0.3517  & 0.01252  & 91.6& 76.4\\
 \hline
 2048 & 1  &  312 & 0.2824& 0.00175  & 93.2& 91.6  &  190 & 0.3247& 0.00299  & 92.5& 91.1 \\
      & 2   & 148&  0.2856 & 0.00151 & 92.6& 90.5& 90 &  0.3256 & 0.00275    & 91.6&88.6\\
      & 4   & 72&  0.2847  & 0.00148  & 93.1& 89.1& 44&0.3272 & 0.00290 & 91.1&85.0 \\
      & 8  & 35& 0.2852  & 0.00153  & 93.2& 87.1 & 21& 0.3283  & 0.00346  & 89.9&80.2\\
      & 16  & 29& 0.2864  & 0.00186  & 93.9& 83.9  & 10& 0.3349  & 0.00458  & 91.2&73.2\\
  \hline
8192 & 1  &  947 & 0.2905& 0.00053  & 93.1& 92.7  &  577& 0.3135& 0.00091  & 93.1& 91.6 \\
      & 2   & 451&  0.2919 & 0.00044 & 93.4& 92.4& 275 &  0.3146 & 0.00084    & 90.7&89.3\\
      & 4   & 219&  0.2925  & 0.00041  & 92.9& 91.5&134&0.3158& 0.00082& 90.7&88.3 \\
      & 8  & 108& 0.2932  & 0.00041  & 93.6& 91.2 &66& 0.3155  & 0.00089  & 90.7&86.0\\
      & 16  & 53& 0.2939  & 0.00043  & 93.0& 88.9  & 32& 0.3184  & 0.00106  & 89.9&82.2\\
 \hline
 \end{tabular}}
\label{tab:result2}
\end{table}


The optimal bandwidth depends on the parameters of the model, thus  it is not possible to obtain this quantity in
practical situations. Indeed, it is well known that  the semiparametric estimators are bandwidth driven estimation
procedures. Due to these peculiarities, an empirical investigation was considered for the bandwidth $m_n =
(\frac{N}{g})^{0.5}$, and the results are in \ref{tab:result3}. It is not surprising that the use of this
bandwidth produces estimates with larger mean-square errors. However, there is an empirical evidence that the reduction
of the mean-square error can be obtained even for  $\nbblock>4$. In addition, the $cr_{r}$ is   more accurate than the
previous case, which can be explained by the fact that the reduction in the number  of  the periodogram
ordinates in the  regression mitigates the effect of the AR coefficient.  The reduction of the impact of the AR coefficient for this
choice  of the bandwidth is also manifested  by the similarity
between the estimates of  the two models.


\begin{table}[h]
\centering \caption{Estimation in ARFIMA($1,d,0$) model using $m_n = (\frac{N}{g})^{0.5}$ } {\footnotesize
\begin{tabular}{c c c c c cc||ccccc }
  \hline
  \hline
    &   &\multicolumn{3}{c}{$d=0.3$, $\phi= -0.3$}&&&\multicolumn{3}{c}{$d=0.3$, $\phi =0.3$}& & \\
    \cline{2-12}
  $N$ & $g$& $m_n$   & $mean$ & $mse$ & $cr_{r}$ &  $cr_{a}$ &$m_n$& $mean$ & $mse$ & $cr_{r}$ & $cr_{a}$\\
  \hline
  512 & 1  &  22 & 0.3016 & 0.02873  & 95.0&88.7  &  22& 0.3116 & 0.02782  & 96.0& 90.0\\
      & 2   & 16 &  0.3039  & 0.01656 & 95.8& 87.1 & 16 &  0.3205 & 0.01781   & 94.6& 86.2\\
      & 4   & 11& 0.2916  & 0.01388 & 94.5&  81.8 & 11& 0.3336  & 0.01441 & 93.6& 80.3\\
   \hline
 2048 & 1  &  45 & 0.3008& 0.01197  & 95.8& 91.4  &  45 & 0.3053& 0.01185  & 94.8& 92.2 \\
      & 2   & 32&  0.3075 & 0.00714 & 95.0& 90.6& 32&  0.3064 & 0.00748    & 94.4&89.2\\
      & 4   & 22&  0.3064  & 0.00470  & 95.3& 88.4& 22&0.3117 & 0.00495 & 95.6&86.9 \\
      & 8  & 16& 0.2999  & 0.00337  & 96.0& 87.2 & 16& 0.3216 & 0.00414 & 93.2&82.8\\
      & 16  & 11&  0.2926  & 0.00305 & 94.7& 83.9  & 11& 0.3390 & 0.00455  & 88.9&72.6\\
  \hline
   8192& 1  &  90  & 0.3066& 0.00573 & 94.7& 91.9&  90  & 0.3076& 0.00536  & 95.0&93.5\\
       & 2   & 64 &  0.3057 & 0.00316 & 94.8& 91.7&  64&  0.3014 & 0.00312    & 95.4& 92.5\\
       & 4   &45& 0.3051  & 0.00217  & 94.4& 90.6  & 45& 0.3057  & 0.00212  & 95.4&91.2\\
       & 8  &32& 0.3046  & 0.00158  & 94.3& 90.1 & 32& 0.3061  & 0.00158  & 94.0&89.5 \\
      & 16  & 22&  0.3045  & 0.00122  & 94.3& 85.0  & 22& 0.3141 & 0.00140  & 92.4&83.5\\
   \hline
\end{tabular}}
\label{tab:result3}
\end{table}

\pagebreak

\section*{Acknowledgements}

The  author V.A. Reisen gratefully acknowledges  partial financial
support from CNPq/ Brazil. G.C. Franco was partially supported by
CNPq-Brazil, and also by {\it Fundação de Amparo à Pesquisa no
Estado de Minas Gerais\/} (FAPEMIG Foundation). The authors also
thank  Giovanni Comarela, undergraduate student under supervision
of V. A. Reisen, who provided the simulation presented in the
paper.

\section{Proof of Theorem \ref{theo:covariance-DFT-multiple-epochs}}
\label{proof:theo:covariance-DFT-multiple-epochs}
Denote by
\begin{equation}
\label{eq:dirichlet-kernel}
D_n(\omega) = \rme^{\rmi \omega} (1 - \rme^{\rmi n \omega})/ (1 - \rme^{\rmi \omega})
\end{equation}
the Dirichlet kernel.
Straightforward calculations show that
\begin{equation}
\label{eq:expression-covariance}
\PE[ d_{n,0}(\omega_j) \bar{d}_{n,\ell}(\omega_k)] =  \int_{-\pi}^{\pi} \dsp(\omega) E_{n,j,k}(\omega) \rme^{-\rmi \ell n \omega} \rmd \omega \eqsp,
\end{equation}
where
\begin{equation}
\label{eq:definition-Enjk}
E_{n,j,k}(\omega) \eqdef (2 \pi n)^{-1} D_n(\omega-\omega_j) \bar{D}_n(\omega-\omega_k) \eqsp.
\end{equation}
The proof then follows from the change of variable $\omega \to n \omega$, using the dominated convergence
theorem.

\section{Proof of Theorem \ref{theo:covariance-DFT-multiple-epochs-bounds}}
\label{sec:theo:covariance-DFT-multiple-epochs-bounds}
We preface the proof by two technical lemmas, which are used in the sequel.
Throughout this section, $C$ is a constant, depending only on  $\mu,\nu$, $\delta_-$ and $\delta_+$, but which may take different values upon each appearance.
\begin{lem}
\label{lem:bound-difference-functions}
There exists a constant $C$ (depending only on $\delta_-$, $\delta_+$, $\mu$) such that, for any $d \in [-\delta_-,\delta_+]$,
$\SMdsp \in \LocalLipshitz{\mu,\nu}$, and $\omega,\omega' \in [-\pi,\pi] \setminus \{0\}$,
\begin{equation}
\label{eq:bound-difference-functions}
\left| |\omega|^{-2d} \SMdsp(\omega) - |\omega'|^{-2d} \SMdsp(\omega') \right| \leq C \SMdsp(0) (|\omega| \wedge |\omega'|)^{-1-2d} |\omega - \omega'| \eqsp.
\end{equation}
In addition, there exists a constant $C$, depending only on $\mu$ and $\nu$ such that, for all integers $\ell \geq 1$ and $j \in \{1, \dots, \tilde{n}/2\}$,
\begin{multline}
\label{eq:bound-difference-integrals}
\int_{-\omega_j}^{\omega_j} |\dsp(\omega) - \dsp(\omega+\pi/\ell n)| \rmd \omega   \\
\leq C \SMdsp(0) \left[ (\ell n)^{-1} \left\{ (\ell n)^{2d} + (j/n)^{-2d} \right\} + (\ell n)^{-\nu} (j/n)^{1-2d} \right] \eqsp.
\end{multline}
\end{lem}

\begin{proof}
The proof of \eqref{eq:bound-difference-functions} is obvious and is omitted for brevity.
Note first that, for any $\omega,\omega' \in [-\pi,\pi]$,
\begin{equation}
\label{eq:bound-difference-dsp}
|\dsp(\omega) - \dsp(\omega')| \leq C \SMdsp(0) \left\{ \left| |\omega|^{-2d} - |\omega'|^{-2d} \right| + |\omega|^{-2d}
\left| \, |\omega| - |\omega'| \, \right|^{\nu}  \right\} \eqsp.
\end{equation}
Applying this inequality with $\omega'= \omega + \pi/\ell n$, yields to
\begin{align*}
&\int_{-\omega_j}^{\omega_j}  |\dsp(\omega) - \dsp(\omega+\pi/\ell n)| \rmd \omega \\
& \leq C \SMdsp(0)
\left\{ \int_{-\omega_j}^{\omega_j} \left| |\omega|^{-2d} - |\omega+ \pi / \ell n|^{-2d} \right| \rmd \omega +
(\ell n)^{-\nu} \int_{-\omega_j}^{\omega_j} |\omega|^{-2d} \rmd \omega \right\} \nonumber\\
& = C \SMdsp(0) \left\{
 \int_{-\omega_j}^{\omega_j} \left| |\omega|^{-2d} - |\omega+ \pi / \ell n|^{-2d} \right| \rmd \omega +
(1-2d)^{-1} (\ell n)^{-\nu} |\omega_j|^{1-2d} \right\}  \eqsp . \nonumber
\end{align*}
On the
intervals $[-\omega_j,-2 \pi/ \ell n]$ and $[\pi/ \ell n, \omega_j]$, we use the bound
\[
\left| |\omega|^{-2d} - |\omega+\pi/\ell n|^{-2d} \right| \leq C |d| \left\{ |\omega|^{-1-2d} + |\omega+\pi/\ell n|^{-1-2d} \right\}  (\ell n)^{-1} \eqsp,
\]
which yields
\[
\int_{- \omega_j}^{-2 \pi/\ell n} + \int_{\pi/ \ell n}^{\omega_j} \left| |\omega|^{-2d} - |\omega+\pi/\ell n|^{-2d} \right| \rmd \omega \leq C \left\{ (\ell n)^{2d} +  n^{2d} j^{-2d} \right\} (\ell n)^{-1} \eqsp.
\]
On the interval $[-2 \pi/ \ell n, \pi/ \ell n]$, we use the bound
$$
\left| |\omega|^{-2d} - |\omega+ \pi/\ell n)|^{-2d} \right| \leq C \left\{ |\omega|^{-2d} + | \omega + \pi/ \ell n|^{-2d} \right\} \eqsp,
$$
which yields
\[
\int_{-\pi/2 \ell n}^{\pi/\ell n} \left| |\omega|^{-2d} - |\omega+ \pi/\ell n)|^{-2d} \right| \rmd \omega \leq C (1-2d)^{-1} (\ell n)^{-1 + 2d} \eqsp,
\]
which concludes the proof.
\end{proof}
Define
\begin{equation}
\label{eq:definition-Delta}
\Delta_{\ell,n}(\omega) \eqdef (2 \pi n)^{-1} \left\{ D_n(\omega + \pi/ \ell n) - D_n(\omega) \right\} \eqsp,
\end{equation}
where $D_n$ is the Dirichlet kernel defined in \eqref{eq:dirichlet-kernel}.
\begin{lem}
\label{lem:modulus-continuity-dirichlet}
%
There exists a constant $C$ such that, for all $n,\ell\geq1$ and $\omega$ such that $0 < \omega \leq \omega+\pi/(\ell n) \leq \pi$,
\begin{equation}
\label{eq:bound-delta}
\left| \Delta_{\ell,n}(\omega) \right| \leq  C \ell^{-1} (1 + n |\omega|)^{-1} 
 \eqsp.
\end{equation}
\end{lem}
\begin{proof}
For any $\omega, \omega' \in(0,\pi]$,
\begin{align*}
|D_n(\omega) - D_n(\omega')| & = \left| \sum_{k=1}^n \left\{ \rme^{\rmi \omega k} - \rme^{\rmi \omega' k} \right\}\right| = \left|  \int_{\omega}^{\omega'} \sum_{k=1}^n k  \; \rme^{\rmi k \lambda} \rmd \lambda \right| \\
& \leq   \int_{\omega}^{\omega'} \left| \sum_{k=1}^n k  \; \rme^{\rmi k \lambda}  \right| \rmd \lambda
\leq C \int_\omega^{\omega'} \frac{n^2} {1+n\lambda} \, \rmd \lambda \\
& = C n\left| \log(1 + n \omega') - \log(1+ n \omega) \right| \eqsp .
\end{align*}
Thus, if $0 < \omega \leq \omega+\pi/(\ell n) \leq \pi$,
\begin{align*}
|\Delta_{\ell,n}(\omega)|
& \leq C \left| \log(1 + n \omega + \pi/\ell) - \log(1+ n \omega) \right| \leq C' (1+n\omega)^{-1} \ell^{-1} \eqsp .
\end{align*}

\end{proof}

\begin{proof}[Proof of Theorem \ref{theo:covariance-DFT-multiple-epochs-bounds}]
The proof for $\ell=0$ follows from \cite{moulines:soulier:1999} (which is a refinement of \cite[Theorem 2]{robinson:1995}). We consider  the case $\ell \ne 0$. Since the function
$\omega \mapsto \dsp(\omega) E_{n,j,k}(\omega)  \rme^{-\rmi \ell n \omega}$ is $2 \pi$-periodic
and $\rme^{-\rmi \ell n (\omega + \pi/ \ell n)}= - \rme^{- \rmi \ell n \omega}$, we may rewrite \eqref{eq:expression-covariance} as
\begin{align*}
&2 \PE[ d_{n,0}(\omega_j) \bar{d}_{n,\ell}(\omega_k)] \\
&\quad =  \int_{-\pi}^{\pi} \left\{ \dsp(\omega) E_{n,j,k}(\omega) - \dsp(\omega + \pi/\ell n) E_{n,j,k}(\omega+\pi/ \ell n) \right\}  \rme^{-\rmi \ell n \omega} \rmd \omega \\
& \quad = A(n,j,k) + B(n,j,k) \eqsp,
\end{align*}
where $E_{n,j,k}(\omega)$ is defined in \eqref{eq:definition-Enjk} and the two terms $A(n,j,k)$ and $B(n,j,k)$ are respectively defined by
\begin{align}
\label{eq:definition-A}
& A(n,j,k) \eqdef  \int_{-\pi}^{\pi} \left\{ \dsp(\omega)  - \dsp(\omega + \pi/ \ell n) \right\} E_{n,j,k}(\omega)   \rme^{-\rmi \ell n \omega} \rmd \omega \eqsp, \\
\label{eq:definition-B}
& B(n,j,k) \eqdef  \int_{-\pi}^{\pi} \dsp(\omega + \pi/ \ell n)   \left\{ E_{n,j,k}(\omega) - E_{n,j,k}(\omega + \pi / \ell n) \right\}  \rme^{-\rmi \ell n \omega} \rmd \omega \eqsp.
\end{align}
First consider $A_{n,jk}$. Denote $\SMdspell = \SMdsp(w)-\SMdsp(\omega+\pi(\ell n))$.

We proceed like in the proof of \cite[Theorem 2]{robinson:1995}. We  decompose $A$ as the following sum
\begin{align*}
A(n,j,k) &\eqdef \sum_{i=1}^6 \int_{W_i(n,j,k)} \left\{ \dsp(\omega)  - \dsp(\omega + \pi/\ell n) \right\} E_{n,j,k}(\omega)   \rme^{-\rmi \ell n \omega} \rmd \omega \\
&\eqdef \sum_{i=1}^6 q_i(n,j,k) \eqsp,
\end{align*}
where $W_1(n,j,k) \eqdef \{ -\omega_j/2 \leq \omega \leq \omega_j/2 \}$, $W_2(n,j,k) \eqdef \{ \omega_j/2 \leq \omega \leq (\omega_j+\omega_k)/2\}$, $W_3(n,j,k) \eqdef \{ (\omega_j+\omega_k)/2 \leq \omega \leq 3 \omega_k/2 \}$,
$W_4(n,j,k) \eqdef \{ 3 \omega_k/2 \leq \omega \leq \pi \}$,  $W_5(n,j,k) \eqdef \{ -\pi \leq \omega \leq -\omega_k\}$, and $W_6(n,j,k) \eqdef \{ - \omega_k \leq \omega \leq - \omega_j/2 \}$.
Note that, for $\omega_0 < \pi$, there exists a constant $C < \infty$ (depending only on $\omega_0$)
such that, for all $\omega \in [-\omega_0,\omega_0]$,
\begin{align}
\label{eq:bound-dirichlet}
&| D_n(\omega) | \leq C n (1 + n |\omega|)^{-1} \eqsp, \\
\label{eq:bound-double-dirichlet}
&| E_{n,j,k}(\omega) | \leq C  n (1 + n |\omega-\omega_j|)^{-1} (1 + n |\omega-\omega_k|)^{-1}.
\end{align}
For $\omega \in W_1(n,j,k)$, \eqref{eq:bound-double-dirichlet} implies that
$|E_{n,j,k}(\omega)| \leq C n k^{-1} j^{-1}$. Using the bound \eqref{eq:bound-difference-integrals}, 
we therefore obtain
\begin{align*}
\omega_j^d \omega_k^d |q_1(n,j,k)| 1 &\leq  j^{d-1} k^{d-1} n^{1-2d} \int_{-\omega_j/2}^{\omega_j/2} \left| \dsp(\omega) - \dsp(\omega+  \pi/ \ell n) \right| \rmd \omega  \nonumber \\
&  \leq  C \SMdsp(0) \,\left( \ell^{-(1-2d)} j^{d-1} k^{d-1} + \ell^{-1}  j^{-d-1} k^{d-1}
+ (\ell n)^{-\nu} j^{-d} k^{d-1} \right) \eqsp.
\end{align*}
For $d \in [0,1/2)$, we have $j^{d-1} \leq j^{-d}$; for $d  \in [0,1/2)$, $j^{-d} k^{d-1} \leq j^{-|d|} k^{|d|-1}$.  
Therefore, for $d \in (-1/2,1/2)$,
\begin{equation}
\label{eq:bound-q1}
\omega_j^d \omega_k^d |q_1(n,j,k)|  \leq
C \SMdsp(0) \, \left\{ \ell^{-\{(1-2d) \wedge 1\}} + (\ell n)^{-\nu} \right\} \; j^{-|d|} k^{|d|-1} \eqsp.
\end{equation}

By \eqref{eq:bound-difference-dsp}, for $\omega \in [\omega_j/2, (\omega_j+\omega_k)/2]$,
$$
|\dsp(\omega) - \dsp(\omega+\pi/\ell n)|  \leq C \SMdsp(0) n^{2d} j^{-1-2d} \ell^{-1} \eqsp.
$$
For  $\omega \in W_2(n,j,k)$, we use the bounds
$|D_n(\omega-\omega_k)| \leq C n \{1+(k-j)\}^{-1}$, and
$$
\int_{\omega_j/2}^{(\omega_k+\omega_j)/2} |D_n(\omega-\omega_j)| \rmd \omega \leq C \log(k)
$$
(see \cite[Lemma 5]{robinson:1994}). Thus
\[
\omega_j^d \omega_k^d |q_2(n,j,k)| \leq C \SMdsp(0) \, \ell^{-1} j^{-1-d} k^d  \{1+(k-j)\}^{-1} \log(k) \eqsp.
\]
We consider separately the cases $j \leq k \leq 2j$ ($j$ and $k$ are
close) and $j > 2k$ ($j$ and $k$ are far apart). If $j \leq k \leq
2j$, then $j^{-1} \leq 2 k^{-1}$, so that $j^{-1-d} k^d \{1 +
(k-j)\}^{-1} \leq 2 j^{-d} k^{d-1}$. On the other hand, if $k \geq
2j$, $j \leq k/2$, $(k-j)^{-1} \leq 2 k^{-1}$. Therefore, $j^{-1-d}
k^d \{1+(k-j)\}^{-1} \leq 2 j^{-d} k^{d-1} $. Combining these two
inequalities and using, for $d < 0$ that $j^{-d} k^{d-1} \leq j^{d}
k^{-d-1}$ yields, for $d \in (-1/2,1/2)$,
\begin{equation}
\label{eq:bound-q2}
\omega_j^d \omega_k^d |q_2(n,j,k)| \leq C \SMdsp(0) \, \ell^{-1} j^{-|d|} k^{|d|-1} \log(k) \eqsp.
\end{equation}
The bound for $q_3$ can be obtained exactly along the same lines. For $\omega \in W_4(n,j,k)$, we use  \eqref{eq:bound-difference-functions}
which shows that
\[
\left|\dsp(\omega) - \dsp(\omega+\pi/\ell n) \right|  \leq C \SMdsp(0) \ell^{-1} n^{2d} k^{-1-2d}  \eqsp,
\]
and, by \eqref{eq:bound-double-dirichlet}, $|E_{n,j,k}(\omega)| \leq C n^{-1} (\omega - \omega_k)^{-2}$, which imply
\begin{multline*}
\omega_j^d \omega_k^d |q_4(n,j,k) | \leq C \SMdsp(0) \ell^{-1}  j^d k^{-1-d} n^{-1} \int_{3 \omega_k/2}^\pi (\omega-\omega_k)^{-2} \rmd \omega  \\
\leq C \SMdsp(0) \ell^{-1} j^d k^{-2-d} \eqsp.
\end{multline*}
For $d \geq 0$, $j^d k^{-2-d} \leq j^{-d} k^{d-1}$ and for $d < 0$, $j^d k^{-2-d} \leq j^d k^{-d-1}$. Therefore,
for any $d \in (-1/2,1/2)$,
\begin{equation}
\label{eq:bound-q4}
\omega_j^d \omega_k^d |q_4(n,j,k)| \leq C \SMdsp(0) \, \ell^{-1} j^{-|d|} k^{|d|-1} \log(k) \eqsp.
\end{equation}
The bound for $q_5$ can be obtained exactly along the same lines.
For $\omega \in W_6(n,j,k)$, we use the bounds $|\dsp(\omega) - \dsp(\omega+\pi/\ell n)|\leq C \SMdsp(0) \ell^{-1} n^{2d} j^{-1-2d}$, $|D_{n}(\omega-\omega_k)| \leq C n k^{-1}$, and
$\int_{-\omega_k}^{-\omega_j/2} |D_n(\omega-\omega_j)| \rmd \omega \leq C \log(k)$, which imply
\[
\omega_j^d \omega_k^d |q_6(n,j,k) | \leq C \SMdsp(0) \, \ell^{-1} j^{-1-d} k^{d-1} \log(k) \eqsp.
\]
For $d \geq 0$, $j^{-1-d} \leq j^{-d}$ and for $d \in (-1/2,0)$, $j^{-1-d} k^{d-1} \leq j^d k^{-d-1}$, which implies,
for $d \in (-1/2,1/2)$ that
\begin{equation}
\label{eq:bound-q6}
\omega_j^d \omega_k^d |q_6(n,j,k)| \leq C \SMdsp(0) \, \ell^{-1} j^{-|d|} k^{|d|-1} \eqsp.
\end{equation}
By combining the bounds obtained in \eqref{eq:bound-q1}, \eqref{eq:bound-q2}, \eqref{eq:bound-q4} and \eqref{eq:bound-q6},
we therefore obtain the following bound
\begin{equation}
\label{eq:bound-A}
\omega_j^d \omega_k^d |A(n,j,k)| \leq C \SMdsp(0) \left\{ \ell^{-(1-2d) \wedge 1} + (\ell n)^{-\nu}  + \ell^{-1} \log(k) \right\} j^{-|d|} k^{|d|-1} \eqsp.
\end{equation}
We now consider the second term $B(n,j,k)$ defined in \eqref{eq:definition-B}. Note that
\begin{multline}
\label{eq:decomposition-Enjk}
E_{n,j,k}(\omega) - E_{n,j,k}(\omega + \pi / \ell n) \\
=  \Delta_{\ell,n}(\omega-\omega_j) \bar{D}_n(\omega - \omega_k) + D_n(\omega + \pi/ \ell n - \omega_j) \bar{\Delta}_{\ell,n}(\omega-\omega_k) \eqsp,
\end{multline}
where $D_n$ is the Dirichlet kernel and $\Delta_{\ell,n}$ is defined in \eqref{eq:definition-Delta}.
In addition, since $D_n(\omega)$ and $\Delta_{\ell,n}(\omega)$ are polynomial in $\rme^{\rmi \omega}$ of degree at most $n$, for any $\ell > 0$, and any $1 \leq j, k \leq [n/2]$,
\begin{equation}
\label{eq:orthogonality-kernel}
\int_{-\pi}^{\pi} \Delta_{\ell,n}(\omega-\omega_j) \bar{D}_n(\omega - \omega_k) \rme^{-\rmi \ell n \omega} \rmd \omega= 0 \eqsp.
\end{equation}
Using this identity together with \eqref{eq:decomposition-Enjk}, we may split $B(n,j,k)$ into two terms $B_1(n,j,k)$ and $B_2(n,j,k)$
which are defined as follows:
\begin{align*}
&B_1 \eqdef  \int_{-\pi}^{\pi} \left\{ \dsp\left(\omega + \frac{\pi}{\ell n}\right) - \dsp\left(\omega_j + \frac{\pi}{\ell n}\right) \right\} \Delta_{\ell,n}(\omega-\omega_j) \bar{D}_n(\omega - \omega_k) \rme^{- \rmi \ell n \omega} \rmd \omega\eqsp, \\
&B_2 \eqdef  \int_{-\pi}^{\pi} \left\{ \dsp\left(\omega + \frac{\pi}{\ell n}\right) - \dsp\left(\omega_k + \frac{\pi}{\ell n}\right) \right\} D_{n}(\omega-\omega_j)  \bar{\Delta}_{\ell,n}(\omega - \omega_k) \rme^{- \rmi \ell n \omega}  \rmd \omega\eqsp.
\end{align*}
These two terms can be handled similarly. We consider $B_1(n,j,k)$. We decompose this term as the sum $B_1(n,j,k) \eqdef \sum_{i=1}^6 \tilde{q}_i(n,j,k)$, with
\begin{multline*}
\tilde{q}_i(n,j,k) = \\
\int_{W_i(n,j,k)} \left\{ \dsp\left(\omega + \frac{\pi}{\ell n}\right) - \dsp\left(\omega_j + \frac{\pi}{\ell n} \right) \right\} \Delta_{\ell,n}(\omega-\omega_j) \bar{D}_n(\omega - \omega_k)  \rme^{-\rmi \ell n \omega} \rmd \omega
\eqsp,
\end{multline*}
where the intervals $W_i(n,j,k)$, $i=1,6$ are defined as above.
For $\omega \in W_1(n,j,k)$, we have
$$\left| \dsp(\omega + \pi/\ell n) - \dsp(\omega_j + \pi/ \ell n) \right| \leq C \SMdsp(0) \, (|\omega|^{-2d} + \omega_j^{-2d}) \eqsp.$$
Moreover, Lemma \ref{lem:modulus-continuity-dirichlet} shows that
$|\Delta_{\ell,n}(\omega-\omega_j)| \leq C \ell^{-1} j^{-1}$ and \eqref{eq:bound-dirichlet} implies  $|D_n(\omega-\omega_k)| \leq C n k^{-1}$.
Therefore,
$$
|\tilde{q}_1(n,j,k)| \leq C \SMdsp(0)  n  (j k \ell)^{-1} \int_{-\omega_j/2}^{\omega_j/2}  (|\omega|^{-2d} + \omega_j^{-2d})  \rmd \omega \leq C \SMdsp(0) n^{2d} j^{-2d} k^{-1} \ell^{-1} \eqsp.
$$
Since for $d < 0$, $j^{-d} k^{d-1} \leq j^d k^{-d-1}$, this bound implies, for any $d \in (-1/2,1/2)$,
\begin{equation}
\label{eq:bound:tildeq-1}
\omega_j^d \omega_k^d |\tilde{q}_1(n,j,k)| \leq C \SMdsp(0) \, \ell^{-1} j^{-|d|} k^{|d|-1} \eqsp.
\end{equation}
Consider now $\omega \in W_2(n,j,k)$. We consider separately the
case $j \leq k \leq 2j$ ($j$ and $k$ are close) and $k > 2j$ ($j$
and $k$ are far apart).  In the first case ($j$ and $k$ close), we
use the bounds 
$$ |\dsp(\omega + \pi /\ell n) - \dsp(\omega_j+\pi/\ell n)| \leq C \SMdsp(0) \omega_j^{-1-2d} |\omega - \omega_j|
$$
(see \eqref{eq:bound-difference-functions}), 
$|\Delta_{\ell,n}(\omega-\omega_k)| \leq C \ell^{-1} n^{-1}
|\omega-\omega_j|^{-1}$ (Lemma \ref{lem:modulus-continuity-dirichlet}), and $|D_n(\omega-\omega_k)| \leq C n\{1
+(k-j)\}^{-1}$, which imply that
\begin{align*}
\omega_j^d \omega_k^d |\tilde{q}_2(n,j,k)|
&\leq C \SMdsp(0) \, \omega_j^d \omega_k^d  \{1 +(k-j)\}^{-1}  \ell^{-1}  \int_{\omega_j/2}^{(\omega_k+\omega_j)/2} \omega_j^{-1-2d} \rmd \omega \\
&\leq C \SMdsp(0) \, \ell^{-1} j^{-1-d} k^d \eqsp.
\end{align*}
Since $j^{-1} \leq 2k^{-1}$, $j^{-1-d} k^{d} \leq 2 j^{-d} k^{d-1} \leq 2 j^{-|d|} k^{|d|-1} $.
Therefore, for $j \leq k \leq 2j$ and all $d \in (-1/2,1/2)$,
\begin{equation}
\label{eq:bound:tildeq-2:close}
\omega_j^d \omega_k^d |\tilde{q}_2(n,j,k)| \leq C \SMdsp(0) \ell^{-1} j^{-|d|} k^{|d|-1} \eqsp.
\end{equation}
In the second case ($j$ and $k$ far apart), we use the bounds 
$
|\dsp(\omega + \pi /\ell n) - \dsp(\omega_j+\pi/\ell n)| \leq C \SMdsp(0) (\omega_k^{-2d} + \omega_j^{-2d})
$, 
$|D_n(\omega-\omega_k)| \leq C n k^{-1}$, and
$$
\int_{\omega_j/2}^{(\omega_k+\omega_j)/2} |\Delta_{\ell,n}(\omega-\omega_j)| \rmd \omega   \leq
C \ell^{-1} \int_{\omega_j/2}^{(\omega_k+\omega_j)/2} (1 + n \omega)^{-1} \rmd \omega \leq C \ell^{-1} n^{-1} \log(k) \eqsp.
$$
This implies that
\begin{align}
\label{eq:bound:tildeq-2:far}
\omega_j^d \omega_k^d |\tilde{q}_2(n,j,k)| &\leq  C \SMdsp(0) \ell^{-1} \left(j^{-d} k^{d-1} + j^d k^{-d-1} \right) \log(k)  \nonumber \\
&\leq C \SMdsp(0) \, \ell^{-1} j^{-|d|} k^{|d|-1} \log(k) \eqsp,
\end{align}
For $\omega \in  W_3(n,j,k)$, and $j \leq k \leq 2j$, we use the bounds $|\dsp(\omega + \pi / \ell n) - \dsp(\omega_j+\pi/\ell n)| \leq C \omega_j^{-1-2d} |\omega-\omega_j|$,
$|\Delta_{\ell,n}(\omega-\omega_j)| \leq C \ell^{-1} n^{-1} |\omega-\omega_j|^{-1}$, and
\begin{equation}
\label{eq:borne-dirichlet-integral}
\int_{(\omega_j+\omega_k)/2}^{3 \omega_k/2} |D_n(\omega-\omega_k)| \rmd \omega \leq C \log(k) \eqsp.
\end{equation}
Therefore, since $k \leq 2 j$,
$$
\omega_j^d \omega_k^d |\tilde{q}_3(n,j,k)| \leq C \SMdsp(0) \log(k) \omega_j^{-1-d} \omega_k^d \ell^{-1} n^{-1} \leq  C \SMdsp(0) \ell^{-1} j^{-|d|} k^{|d|-1} \eqsp.
$$
For $\omega \in  W_3(n,j,k)$, and $j < k/2$, we use the bounds $|\dsp(\omega + \pi / \ell n) - \dsp(\omega_j+\pi/\ell n)| \leq C \SMdsp(0) \left( \omega_k^{-2d} + \omega_j^{-2d} \right)$,  $|\Delta_{\ell,n}(\omega-\omega_j)| \leq C \ell^{-1} n^{-1} \omega_k^{-1}$, and \eqref{eq:borne-dirichlet-integral}.
Therefore, for $j < k/2$,
\begin{multline*}
\omega_j^d \omega_k^d |\tilde{q}_3(n,j,k)| \leq C \SMdsp(0) \log(k) \ell^{-1} k^{-1} \left( (k/j)^{d} + (j/k)^{d}  \right) \\ \leq  C \SMdsp(0) \ell^{-1} j^{-|d|} k^{|d|-1} \eqsp.
\end{multline*}
For $\omega \in W_4(n,j,k)$, we use the bounds $|\dsp(\omega + \pi /\ell n) - \dsp(\omega_j+\pi/\ell n)| \leq C \SMdsp(0) (\omega^{-2d} + \omega_j^{-2d})$,
$|\Delta_{\ell,n}(\omega-\omega_j)| \leq C \ell^{-1} n^{-1} \omega^{-1}$, and $|D_n(\omega-\omega_k)| \leq C \omega^{-1}$,
which implies that
\begin{multline*}
|\tilde{q}_4(n,j,k)| \leq C \SMdsp(0)  \ell^{-1} n^{-1} \int_{3\omega_k/2}^\pi (\omega^{-2d}+\omega_j^{-2d})\omega^{-2} \rmd \omega
\\ 
\leq C \SMdsp(0) \ell^{-1} k^{-1} (\omega_k^{-2d}+\omega_j^{-2d}) \eqsp .
\end{multline*}
Hence,
\begin{equation}
\label{eq:bound:tildeq-4}
\omega_j^d \omega_k^d |\tilde{q}_4(n,j,k)| \leq C \SMdsp(0) \ell^{-1} k^{-1} \left( (k/j)^{d} + (j/k)^{d}  \right)
\leq C \SMdsp(0) \ell^{-1} j^{-|d|} k^{|d|-1} \eqsp.
\end{equation}
The bound for $\tilde{q}_5$ can be obtained exactly along the same lines. For $\omega \in W_6(n,j,k)$, we use the bounds $|\dsp(\omega + \pi / \ell n) - \dsp(\omega_j+\pi/\ell n)| \leq C \SMdsp(0) (\omega_k^{-2d} + \omega_j^{-2d})$,
$\int_{W_6(n,j,k)} |\Delta_{\ell,n}(\omega-\omega_j)| \rmd \omega \leq C \ell^{-1} n^{-1} \log(k)$ and $|D_n(\omega-\omega_k)| \leq C n k^{-1}$,  which implies that
\begin{multline*}
\omega_j^d \omega_k^d |\tilde{q}_6(n,j,k)| \leq C \SMdsp(0) \ell^{-1} \left( j^{-d} k^{d-1} + j^{d} k^{-d-1} \right) \log(k)
\\ \leq C \SMdsp(0) \ell^{-1} j^{-|d|} k^{|d|-1} \log(k) \eqsp.
\end{multline*}
This concludes the  proof.
\end{proof}
\section{Proof of Theorem \ref{theo:gphgauss}}
\label{sec:proof:theo:gphgauss}
The proof is similar to the one used for the pooled periodogram given in \cite{moulines:soulier:2003} (which simplifies the arguments given \cite{robinson:1995} being based on the
central limit theorem for non-linear functions of Gaussian vectors given in \cite{soulier:2001}).
The error $\hat{d}_{\nbblock,n}-d$ is naturally decomposed into a stochastic and a bias terms as follows
\begin{multline*}
\hat{d}_{\nbblock,n}-d = \sum_{k=1}^{m_n} a_{k,n}(m_n) \log \left[\bar{I}_{\nbblock,n}(\omega_k)/\dsp(\omega_k)\right] \\
+ \sum_{k=1}^{m_n} a_{k,n}(m_n) \log \left[ \SMdsp(\omega_k) / \SMdsp(0) \right] = S_n(m_n,\nbblock) + B_n(m_n,\SMdsp) \eqsp,
\end{multline*}
where the coefficients $\{ a_{k,n}(m) \}$ is defined in \eqref{eq:definition-weight-GPH}.
In the previous expression, $S_n(m_n,\nbblock)$ is a stochastic fluctuation term (which depends in particular
on the number of epochs) and $B_n(m_n,\SMdsp)$ is the bias caused by the approximation in the neighborhood
of the zero frequency of $\SMdsp$ by a constant. The result will follow from the weak convergence
of the stochastic term $S_n(m_n,\nbblock)$ and from a bound for the bias term $B_n(m_n,\SMdsp)$. By \cite[Lemma 6.1]{moulines:soulier:2003}, there exists a
constant $C(\Omega_0,\beta,\mu)$ such that, for any $\SMdsp \in \LocalHolder{\Omega_0}{\beta}{\mu}$ and any nonnegative integer $m$ satisfying $2 \pi m/n \leq \Omega_0$,
$|B_n(m,\SMdsp)| \leq C(\Omega_0,\beta,\mu) (m/n)^\beta$. Therefore, under the condition \eqref{eq:ratemgph}, $\lim_{n \to \infty} \sqrt{m_n} B_n(m_n,\SMdsp)=0$.

 To establish the weak convergence result, we use
\cite[Theorem 9.7]{moulines:soulier:2003} (which is based on
\cite[Theorem 4.1]{soulier:2001}). To simplify the notations, put
$a_{k,n}= a_{k,n}(m_n)$. Since $\lim_{n \to \infty}
m_n^{-1} + n^{-1} m_n= 0$, it follows from \cite{robinson:1995} that
\[ \lim_{n \to \infty} m_n \sum_{k=1}^{m_n} a_{k,n}^2= 1/4 \quad \text{and} \quad \max_{1 \leq k \leq m_n} |a_{k,n}| = O\left(m_n^{-1/2} \log(m_n) \right) \eqsp. \]
Set $v_n= \lfloor m_n^{\gamma} \rfloor$ for some $\gamma \in  (1/2,1)$ and $\bar{g}_{n} \eqdef m_n^{-1} \sum_{k=1}^{m_n} [-2\log(\omega_k)]$, where the
function $g$ is defined in \eqref{eq:definition-weight-GPH}. Note that
\begin{align*}
&\max_{1 \leq k \leq m_n} |a_{k,n}| \sum_{k=1}^{v_n} |a_{k,n}| \log^2(v_n) \leq m_n^{\gamma}  \log^2(m_n) \left(\max_{1 \leq k \leq m_n} |a_{k,n}|\right)^2 \to 0 \eqsp, \\
&\sum_{k=1}^{m_n} |a_{k,n}| \log(v_n)/v_n \leq m_n^{1/2-\gamma} \log(m_n) \frac{m_n^{-1} \sum_{k=1}^{m_n} \left| -2\log(\omega_k) - \bar{g}_{n} \right|}{m_n^{-1} \sum_{k=1}^{m_n} \left\{ -2\log(\omega_k) - \bar{g}_{n} \right\}^2} \to 0 \eqsp,
\end{align*}
To apply \cite[Theorem 9.7]{moulines:soulier:2003}, we finally need to prove that there exists a constant $C$ such that, for all $k \in \{1, \dots, m_n\}$,
$ \PE\left[ \log^2\left(\bar{I}_{\nbblock,n}(\omega_k)/\dsp(\omega_k)\right) \right] \leq C$. Corollary \ref{cor:decomposition-averaged-periodogram}
shows that this bounds holds for any $k \in \{ K, \dots, m_n\}$, where $K$ is a fixed integer. For the first $K$ frequencies, we use Theorem \ref{theo:covariance-DFT-multiple-epochs} to
show that, for any given $k$ and uniformly in $n$, the minimal eigenvalue of the covariance matrix of
the random vector
\[
[\mathrm{Re}\left\{d_{n,0}(\omega_k)\right\}, \mathrm{Im}\left\{d_{n,0}(\omega_k)\right\}, \dots, \mathrm{Re}\left\{d_{n,\nbblock-1}(\omega_k)\right\},\mathrm{Im}\left\{d_{n,\nbblock-1}(\omega_k)\right\}] \]
is bounded away from zero.

%

\end{document}